\documentclass[letterpaper]{amsart}
\usepackage{amsmath,amsthm,amssymb}       
\usepackage{epsfig,graphicx}

\begin{document}

\renewcommand{\labelenumi}{\theenumi.}

\bibliographystyle{amsplain}



\hfuzz4pt 

\setlength{\floatsep}{15pt plus 12pt}
\setlength{\textfloatsep}{\floatsep}

\setlength{\intextsep}{\floatsep}\def\ldotsplus{\mathinner{\ldotp\ldotp\ldotp\ldotp}}
\setlength{\intextsep}{\floatsep}\def\ldotscomm{\mathinner{\ldotp\ldotp\ldotp,}}
\def\fourdots{\relax\ifmmode\ldotsplus\else$\m@th \ldotsplus\,$\fi}
\def\dotcoms{\relax\ifmmode\ldotscomm\else$\m@th \ldotsplus\,$\fi}

\theoremstyle{plain}
\newtheorem{thm}{Theorem}
\newtheorem{lem}[thm]{Lemma}
\newtheorem{cor}[thm]{Corollary}
\newtheorem{prop}[thm]{Proposition}
\newtheorem{conj}[thm]{Conjecture}

\theoremstyle{definition}
\newtheorem{dfn}[thm]{Definition}
\newtheorem{exmp}[thm]{Example}
\newtheorem{remark}[thm]{Remark}
\newtheorem{fact}[thm]{Fact}



\newcommand{\thmref}[1]{Theorem~\ref{#1}}
\newcommand{\corref}[1]{Corollary~\ref{#1}}
\newcommand{\propref}[1]{Proposition~\ref{#1}}
\newcommand{\lemref}[1]{Lemma~\ref{#1}}
\newcommand{\lemsref}[1]{Lemmas~\ref{#1}}
\newcommand{\factref}[1]{Fact~\ref{#1}}
\newcommand{\rmkref}[1]{Remark~\ref{#1}}
\newcommand{\remarkref}[1]{Remark~\ref{#1}}
\newcommand{\defref}[1]{Definition~\ref{#1}}
\newcommand{\egref}[1]{Example~\ref{#1}}
\newcommand{\figref}[1]{Figure~\ref{#1}}
\newcommand{\figsref}[2]{Figures~\ref{#1},\ref{#2}}
\newcommand{\secref}[1]{Section~\ref{#1}}

\newcommand{\miniprefig}
      {\scalebox{.25}{\includegraphics{hexpics/ministr.eps}}}

\newcommand{\minimslfig}
      {\scalebox{.08}{\includegraphics{mslpics/sec3412.eps}}}

\newcommand{\miniinvfig}
      {\scalebox{.08}{\includegraphics{mslpics/sec4231.eps}}}

\newcommand{\minihexafig}
      {\scalebox{.07}{\includegraphics{mslpics/minihexa.eps}}}
\newcommand{\minihexbfig}
      {\scalebox{.07}{\includegraphics{mslpics/minihexb.eps}}}
\newcommand{\minihexcfig}
      {\scalebox{.07}{\includegraphics{mslpics/minihexc.eps}}}
\newcommand{\minihexdfig}
      {\scalebox{.07}{\includegraphics{mslpics/minihexd.eps}}}
\newcommand{\minihexefig}
      {\scalebox{.07}{\includegraphics{mslpics/minihexe.eps}}}
\newcommand{\minihexffig}
      {\scalebox{.07}{\includegraphics{mslpics/minihexf.eps}}}
\newcommand{\minihexgfig}
      {\scalebox{.07}{\includegraphics{mslpics/minihexg.eps}}}

\newcommand{\mymslfig}[2]{\begin{figure}[htbp]\begin{center}
      {\scalebox{.4}{\includegraphics{#1.eps}}}
      \caption{#2}\label{fig:#1}
    \end{center}\end{figure}}

\newcommand{\myintromslfig}[2]{\begin{figure}[htbp]\begin{center}
      {\scalebox{.4}{\includegraphics{#1.eps}}}
      \footnotesize{\textbf{Figure 0.1.} #2}
    \end{center}\end{figure}}

\newcommand{\mysmmslfig}[2]{\begin{figure}[ht]\begin{center}
      {\scalebox{.32}{\includegraphics{#1.eps}}}
      \caption{#2}\label{fig:#1}
    \end{center}\end{figure}}

\newcommand{\myhexfig}[2]{\begin{figure}[ht]\begin{center}
    \input{hexpics/#1.pstex_t}\caption{#2}\label{fig:#1}
    \end{center}\end{figure}}

\newcommand{\mynewhexfig}[2]{\begin{figure}[ht]\begin{center}
    \input{hexpics/#1.pstex_t}\caption{#2}\label{fig:#1}
    \end{center}\end{figure}}

\newcommand{\mynewfig}[2]{\begin{figure}[ht]\begin{center}
      {\scalebox{.3}{\includegraphics{newpics/#1.eps}}}
      \caption{#2}\label{fig:#1}
    \end{center}\end{figure}}

\newcommand{\sfrac}[2]{\genfrac{\{}{\}}{0pt}{}{#1}{#2}}
\newcommand{\sumsb}[1]{\sum_{\substack{#1}}}  
\newcommand{\minsb}[1]{\substack{#1}}  
        
\newcommand{\fsn}{S_n}
\newcommand{\fsp}[1]{S_{#1}}
\newcommand{\frg}{\mathfrak{g}}
\newcommand{\frh}{\mathfrak{h}}
\newcommand{\slgn}{\mathfrak{sl}_n}

\newcommand{\io}{{i_1}}
\newcommand{\iw}{{i_2}}
\newcommand{\inn}{{i_r}}

\newcommand{\bbA}{\mathbb{A}}
\newcommand{\bbR}{\mathbb{R}}
\newcommand{\bbC}{\mathbb{C}}
\newcommand{\bbZ}{\mathbb{Z}}
\newcommand{\bbN}{\mathbb{N}}
\newcommand{\bbQ}{\mathbb{Q}}

\newcommand{\wbi}{w^{\hat{i}}}
\newcommand{\xbi}{x^{\hat{i}}}
\newcommand{\zbi}{z^{\hat{i}}}
\newcommand{\wbio}{w^{\hat{1}}}
\newcommand{\wsbio}{ws^{\hat{1}}}
\newcommand{\xbio}{x^{\hat{1}}}
\newcommand{\zbio}{z^{\hat{1}}}

\newcommand{\tti}{t}
\newcommand{\wti}{{\widetilde{w}}}
\newcommand{\xti}{\protect{\widetilde{x}}}
\newcommand{\uti}{{\widetilde{u}}}
\newcommand{\vti}{{\widetilde{v}}}
\newcommand{\yti}{{\widetilde{y}}}
\newcommand{\zti}{{\widetilde{z}}}

\newcommand{\nts}{\negthickspace}

\newcommand{\trans}{\mathcal{T}}
\newcommand{\simref}{\mathcal{S}}
\newcommand{\calh}{\mathcal{H}}
\newcommand{\calw}{\mathcal{W}}


\newcommand{\msp}{\textsc{msp}}

\newcommand{\br}{\mathbf{a}}
\newcommand{\brp}{\mathbf{a'}}
\newcommand{\bsig}{{\boldsymbol{\sigma}}}
\newcommand{\bnu}{{\boldsymbol{\nu}}}
\newcommand{\bgam}{{\boldsymbol{\gamma}}}
\newcommand{\bdel}{{\boldsymbol{\delta}}}
\newcommand{\bal}{{\boldsymbol{\alpha}}}
\newcommand{\bbe}{{\boldsymbol{\beta}}}

\newcommand{\rank}{\operatorname{lvl}}
\newcommand{\pt}{\operatorname{pt}}
\newcommand{\heap}{\operatorname{Heap}}
\newcommand{\lcz}{\operatorname{lcz}}
\newcommand{\rcz}{\operatorname{rcz}}
\newcommand{\mcz}{\operatorname{mcz}}
\newcommand{\ver}{\operatorname{ver}}
\newcommand{\codim}{\operatorname{codim}}
\newcommand{\bcone}{\operatorname{Cone_{\wedge}}}
\newcommand{\ucone}{\operatorname{Cone^{\vee}}}
\newcommand{\charp}{\operatorname{char}}

\newcommand{\IH}{\operatorname{IH}}

\newcommand{\defeq}{:=}

\newcommand{\digx}{\operatorname{mat}(x)}
\newcommand{\digxt}{\operatorname{mat}(\widetilde{x})}
\newcommand{\digwt}{\operatorname{mat}(\widetilde{w})}
\newcommand{\augdigx}{\operatorname{mat}'(x)}
\newcommand{\digw}{\operatorname{mat}(w)}
\newcommand{\digp}[1]{\operatorname{mat}(#1)}
\newcommand{\dpw}[1]{d_{#1,w}}
\newcommand{\dpb}[2]{d_{#1,#2}}
\newcommand{\duv}{d_{u,v}}
\newcommand{\dvw}{d_{v,w}}
\newcommand{\dzw}{d_{z,w}}
\newcommand{\dxw}{d_{x,w}}
\newcommand{\dyw}{d_{y,w}}
\newcommand{\txw}{\Theta_{x,w}}
\newcommand{\dxtwt}{d_{\widetilde{x},\widetilde{w}}}

\newcommand{\rxjw}{\mathcal{R}(x_i,w)}
\newcommand{\rxjmw}{\mathcal{R}(x_{i-1},w)}
\newcommand{\rxow}{\mathcal{R}(x_0,w)}

\newcommand{\area}{\mathcal{A}}

\newcommand{\ph}{\phi}
\newcommand{\pht}{\phi_t}
\newcommand{\phs}{\phi_s}
\newcommand{\phtp}{\phi_{t'}}
\newcommand{\phtal}{\phi_{t_{\alpha_j}}}

\newcommand{\cP}{{\mathcal P}}
\newcommand{\xsing}{{X_w^{\text{sing}}}}
\newcommand{\maxsing}{\operatorname{\protect{Maxsing}}(X_\protect{w})}
\newcommand{\maxsingt}{\operatorname{Maxsing}(X_{\widetilde{w}})}
\newcommand{\schub}[1]{X_{#1}}

\newcommand{\sij}{s_{i_j}}
\newcommand{\sik}{s_{i_k}}
\newcommand{\sio}{s_{i_1}}
\newcommand{\sir}{s_{i_r}}


\newcommand{\puv}{P_{u,v}}
\newcommand{\pxw}{P_{x,w}}
\newcommand{\pww}{P_{w,w}}
\newcommand{\pyw}{P_{y,w}}
\newcommand{\pyv}{P_{y,v}}
\newcommand{\pzv}{P_{z,v}}
\newcommand{\pxws}{P_{x,ws}}
\newcommand{\pxsws}{P_{xs,ws}}
\newcommand{\pxtwt}{P_{\xti,\wti}}
\newcommand{\pzw}{P_{z,w}}
\newcommand{\pxiwi}{P_{x^{-1},w^{-1}}}
\newcommand{\psxw}{P_{sx,w}}
\newcommand{\pxsw}{P_{xs,w}}
\newcommand{\pxz}{P_{x,z}}
\newcommand{\pkl}[2]{P_{#1,#2}}


\newcommand{\delxw}{\Delta(x,w)}
\newcommand{\delyw}{\Delta(y,w)}
\newcommand{\delytw}{\Delta(yt,w)}

\newcommand{\cpw}{{\mathcal P}(\br)}
\newcommand{\cpww}{{\mathcal P}_w(\br)}
\newcommand{\cpws}{{\mathcal P}(\br s)}
\newcommand{\cpxw}{{\mathcal P}_x(\br)}
\newcommand{\cpyw}{{\mathcal P}_y(\br)}
\newcommand{\cpxsw}{{\mathcal P}_{xs}(\br)}
\newcommand{\cpxsws}{{\mathcal P}_{xs}(\br/s)}
\newcommand{\cpxws}{{\mathcal P}_x(\br/s)}
\newcommand{\cpzero}{{\mathcal P}_x^0(\br)}
\newcommand{\cpone}{{\mathcal P}_x^1(\br)}
\newcommand{\cpeps}{{\mathcal P}_x^\epsilon(\br)}

\newcommand{\cd}{{\mathcal D}}
\newcommand{\Dbr}{\Delta_\bsig}
\newcommand{\Dbrj}{\Delta_{\bsig[j]}}
\newcommand{\Dbrr}{\Delta_{\bsig[r]}}
\newcommand{\Dbrk}{\Delta_{\bsig[k]}}
\newcommand{\Dbrjm}{\Delta_{\bsig[j-1]}}

\newcommand{\apw}{P(\br)}
\newcommand{\apww}{P_w(\br)}
\newcommand{\apxw}{P_x(\br)}
\newcommand{\apxsw}{P_{xs}(\br)}
\newcommand{\apxws}{P_x(\br/s)}
\newcommand{\apxsws}{P_{xs}(\br/s)}

\newcommand{\apzw}{P_z(\br)}
\newcommand{\apzsw}{P_{zs}(\br)}
\newcommand{\apzws}{P_z(\br/s)}
\newcommand{\apzsws}{P_{zs}(\br/s)}

\newcommand{\apew}{P_e(\br)}
\newcommand{\cpew}{{\mathcal P}_e(\br)}


\newcommand{\fl}{\operatorname{fl}}
\newcommand{\ra}{\operatorname{unfl}}
\newcommand{\im}{\operatorname{Im}}
\newcommand{\slide}{\operatorname{\tau}}

\newcommand{\tpq}{t_{p,q}}
\newcommand{\tab}{t_{a,b}}
\newcommand{\talbt}{t_{\alpha,\beta}}
\newcommand{\tac}{t_{a,c}}
\newcommand{\tbc}{t_{b,c}}
\newcommand{\tcd}{t_{c,d}}
\newcommand{\tij}{t_{i,j}}
\newcommand{\tik}{t_{i,k}}
\newcommand{\tjk}{t_{j,k}}
\newcommand{\til}{t_{i,l}}
\newcommand{\talj}{t_{\alpha_j}}
\newcommand{\tgd}{t_{\gamma,\delta}}

\newcommand{\ttpq}{t_{p,q}}
\newcommand{\ttab}{t_{a,b}}
\newcommand{\ttalbt}{t_{\alpha,\beta}}
\newcommand{\ttac}{t_{a,c}}
\newcommand{\ttbc}{t_{b,c}}
\newcommand{\ttcd}{t_{c,d}}
\newcommand{\ttij}{t_{i,j}}
\newcommand{\ttik}{t_{i,k}}
\newcommand{\ttjk}{t_{j,k}}
\newcommand{\ttalj}{t_{\alpha_j}}
\newcommand{\ttgd}{t_{\gamma,\delta}}

\newcommand{\rp}[2]{\mathcal{R}({#1},{#2})}
\newcommand{\qp}[2]{\mathcal{E}_{{#1},{#2}}(x,w)}
\newcommand{\qpt}[2]{\mathcal{E}_{{#1},{#2}}(x,w)}

\newcommand{\rxw}{\mathcal{R}(\protect{x},\protect{w})}
\newcommand{\rytw}{\mathcal{R}(yt,w)}
\newcommand{\rysw}{\mathcal{R}(ys,w)}
\newcommand{\rysiw}{\mathcal{R}(ys_i,w)}
\newcommand{\qabxw}{\mathcal{E}_{a,b}(x,w)}
\newcommand{\qptxw}{\mathcal{E}_t(x,w)}
\newcommand{\qptabxw}{\mathcal{E}_{\tab}(x,w)}
\newcommand{\qbcxw}{\mathcal{E}_{b,c}(x,w)}
\newcommand{\qcdxw}{\mathcal{E}_{c,d}(x,w)}
\newcommand{\qbcyw}{\mathcal{E}_{b,c}(y,w)}

\newcommand{\qabxtwt}{\mathcal{E}_{a,b}(\widetilde{x},\widetilde{w})}
\newcommand{\qptxtwt}{\mathcal{E}_t(\widetilde{x},\widetilde{w})}
\newcommand{\qptabxtwt}{\mathcal{E}_{\tab}(\widetilde{x},\widetilde{w})}
\newcommand{\qbcxtwt}{\mathcal{E}_{b,c}(\widetilde{x},\widetilde{w})}
\newcommand{\qcdxtwt}{\mathcal{E}_{c,d}(\widetilde{x},\widetilde{w})}

\newcommand{\rxiwi}{\mathcal{R}(x^{-1},w^{-1})}
\newcommand{\rysiwi}{\mathcal{R}(sy^{-1},w^{-1})}
\newcommand{\ryw}{\mathcal{R}(y,w)}
\newcommand{\ryiwi}{\mathcal{R}(y^{-1},w^{-1})}
\newcommand{\rxtwt}{\mathcal{R}(\widetilde{x},\widetilde{w})}
\newcommand{\rxtw}{\mathcal{R}(xt,w)}

\newcommand{\ykm}{y_{k,m}}
\newcommand{\xkm}{x_{k,m}}
\newcommand{\xktm}{y_{k,m}}
\newcommand{\xrs}{x_{r,s}}
\newcommand{\xrts}{y_{r,s}}
\newcommand{\wkm}{w_{k,m}}
\newcommand{\vkm}{v_{k,m}}
\newcommand{\wktm}{v_{k,m}}
\newcommand{\wrs}{w_{r,s}}
\newcommand{\wrts}{v_{r,s}}


\title[Inverse Kazhdan-Lusztig polynomials]{Two formulae for inverse
  Kazhdan-Lusztig polynomials in $\boldsymbol{\fsn}$}


\begin{abstract}
  Let $w_0$ denote the permutation $[n,n-1,\ldots,2,1]$.  We give two
  new explicit formulae for the Kazhdan-Lusztig polynomials
  $\pkl{w_0w}{w_0x}$ in $\fsn$ when $x$ is a maximal element in the
  singular locus of the Schubert variety $X_w$.  To do this, we
  utilize a standard identity that relates $\pxw$ and
  $\pkl{w_0w}{w_0x}$.
\end{abstract}

\author{Gregory S. Warrington} \email{warrington@math.umass.edu}

\address{Author's address:
  Dept.\ of Mathematics \& Statistics\\
  University of Massachusetts\\ Amherst, MA 01003}

\maketitle

\section{Introduction}
\label{sec:intro}

Kazhdan-Lusztig (KL) polynomials were introduced by Kazhdan and
Lusztig in \cite{K-L1} in their study of the representations of Hecke
algebras of Coxeter groups.  Since then, these polynomials have been
discovered to have many important interpretations in the context of
Lie theory and Schubert varieties (see \cite{Deodhar94,K-L1,K-L2}).
However, their combinatorial structure is far from clear even though
numerous people have results in specific cases (see \cite{BLak} for
an overview of such results).  In this paper we give explicit formulae
for the KL polynomials $\pxw$ related to certain maximal singular
points in the singular locus of the Schubert variety $X_w$.  These
singular points correspond to the first points where the KL
polynomials are non-trivial.  To state our results precisely, we first
introduce the following four families of permutations:
\begin{dfn}\label{def:permdef}\ \\
  For $k,m \geq 1$, define
  \begin{align*}
    \xkm &= [k,\dotcoms 1,k+m,\dotcoms k+1],\\
    \wkm &= [k+m,k,\dotcoms 2,k+m-1,\dotcoms k+1,1].
  \end{align*}
  For $k,m \geq 1$, define
  \begin{align*}
    \xktm &= [k,\dotcoms 1,k+2, k+1,
    k+2+m,\dotcoms k+2+1],\\
    \begin{split}
      \wktm &= [k+2,k,\dotcoms 2,k+m+2,1,k+2+(m-1),\dotcoms k+2+1,k+1].
    \end{split}
  \end{align*}
\end{dfn}

Writing $w_0$ for the permutation $[n,n-1,\ldots,1]$, our main results
are the following two formulae:
\begin{thm}\ \\\label{thm:mainklinv}
  \begin{enumerate}
  \item For $k,m \geq 1$,
    \begin{equation*}\label{inveq:4231}
      P_{w_0\wkm,w_0\xkm} = \sum_{r=0}^{\min(k-1,m-1)}
      \binom{k-1}{r}\binom{m-1}{r} q^r.
    \end{equation*}
  \item For $k,m \geq 1$,
    \begin{equation*}\label{inveq:3412}
      P_{w_0\wktm,w_0\xktm} = 1 + (k+m-1)q.
    \end{equation*}
  \end{enumerate}
\end{thm}

The pairs $(\xkm,\wkm)$ and $(\xktm,\wktm)$ correspond (using
\thmref{thm:stdfacts}.\ref{fact:pxw}) to two of the three types of
irreducible components of the singular loci of Schubert varieties in
$SL(n)/B$.  While our combinatorial techniques do not easily extend to
the third type, the analogue of \thmref{thm:mainklinv} for this third
type can be found in \cite{Lascoux95}.

\begin{remark}
  Let $u,v$ be elements in an arbitrary Coxeter group and set $c(u,v)$
  to be the number of coatoms in the Bruhat interval $[u,v]$.  Brenti
  shows in \cite{BrentiUpper} that the coefficient of $q$ in $P_{u,v}$
  is bounded above by $c(u,v)-1$.  The intervals
  $[w_0w_{k,k},w_0x_{k,k}]$ afford a class of intervals in $\fsn$ for
  which the coefficient of $q$ in $\pkl{w_0w_{k,k}}{w_0x_{k,k}}$
  asymptotically approaches $c(w_0w_{k,k},w_0x_{k,k})$.  This confirms
  the asymptotic tightness of Brenti's bound.
\end{remark}

\secref{sec:prelim} contains the necessary preliminaries and
\secref{sec:kl-polyn-inverse} contains our proof of
\thmref{thm:mainklinv}.

\section{Preliminaries}
\label{sec:prelim}

We now introduce the necessary background on KL polynomials and
$\fsn$.  The reader is referred to \cite{Hum} for a more leisurely
introduction to most of combinatorial material presented here.  A good
reference for Schubert varieties and the information about them
encoded by KL polynomials is \cite{BLak}.

\subsection{The symmetric group}
\label{sec:symm}

We will view elements of $\fsn$ as permutations on $[1,\dotcoms n]$
with elements $s_i$ of the generating set $\simref =\{s_i\}_{i\in
  [1,\dotcoms n-1]}$ associated with the adjacent transpositions
$(i,i+1)$.  $t_{i,j}$ will denote the transposition $(i,j)$.  We have a
one-line notation for a permutation $w$ given by writing the image of
$[1,\dotcoms n]$ under the action of $w$: $[w(1),w(2),\dotcoms w(n)]$.
The length function for $\fsn$ is given by
\begin{equation*}
  l(w) = |\{1 \leq i < j \leq n: w(i) > w(j)\}|.
\end{equation*}
We will denote the ordered pair of permutations $x$ and $w$ by
$(x,w)$.  Finally, let $w_0 = [n,n-1,\dotcoms 1]\in\fsn$ denote the
element of maximal length in $\fsn$.

We say that $w\in\fsn$ is \textit{$v$-avoiding} for $v\in\fsp{k}$ if we
cannot find $1 \leq i_1 < \cdots < i_k \leq n$ with $w(i_1),\dotcoms
w(i_k)$ in the same relative order as $v(1),\dotcoms v(k)$; i.e., no
submatrix of $\digw$ on rows $i_1,\dotcoms i_k$ and columns
$w(i_1),\dotcoms w(i_k)$ is the permutation matrix of $v$.  There are
many properties pertaining to $\fsn$ and Schubert varieties that can
be characterized efficiently in terms of pattern avoidance (see, e.g.,
\cite{B3,BL,hexagon,LS2,Stem4}).

For $w\in\fsn$ and $1 \leq i_1 < \cdots < i_k \leq n$ for $k\leq n$,
define $\fl[w(i_1),w(i_2),\ldots,w(i_k)]$ to be the unique permutation
$[v(1),\ldots, v(k)] \in \fsp{k}$ such that $v(j) < v(k)$ if and only
if $w(i_j) < w(i_k)$.

We now introduce an important partial order on $\fsn$.  The
characterization we give in \defref{def:bdf} is non-standard, and
requires the following definition, but it is equivalent to more
common descriptions such as the Tableau Criterion.
\begin{dfn}\label{def:rnkdef}
  Let $x,w \in \fsn$, $p,q\in \bbZ$.  Define $r_w(p,q) \defeq |\{i
  \leq p: w(i) \geq q\}|$ and $\dxw(p,q) \defeq r_w(p,q) - r_x(p,q)$.
\end{dfn}

\begin{dfn}\label{def:bdf}
  We define the \textit{Bruhat partial order} ``$\leq$'' on $\fsn$ by
  setting $x\leq w$ if and only if $\dxw(p,q) \geq 0$ for all $p,q$.
\end{dfn}
\begin{lem}\label{lem:trans}
  If $x \leq y \leq w$, then $\dxw - \dyw$ is everywhere non-negative.
\end{lem}

We can now introduce our pictorial version of the Bruhat order.  
\begin{dfn}
  A \textit{Bruhat picture for $x,w\in\fsn$} is an overlay of some of
  the rows and columns of their permutation matrices that is augmented
  by shading the regions on which $\dxw \geq 1$.
\end{dfn}
\figref{fig:brupic} displays examples of this notation.  Let $\digw$
refer to the permutation matrix for $w$.  Entries of $\digx$ (resp.,
$\digw$) are denoted by black disks (resp., open circles).  Positions
corresponding to 1's of both $\digx$ and $\digw$ are denoted by a
black disk and a larger concentric circle.
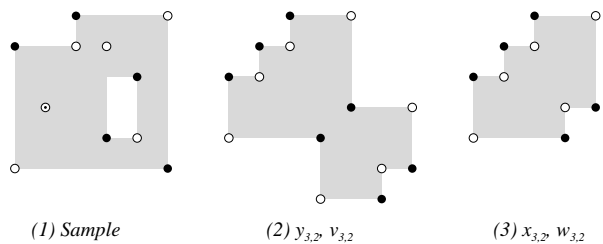
\mysmmslfig{brupic}{1) Bruhat picture for $x = [3,1,5,2,4,6]$, $w =
  [6,3,4,2,5,1]$.  2) and 3) give visualizations of the permutations
  defined in \defref{def:permdef}.}

Our proof of \thmref{thm:mainklinv} is inductive and reduces the
calculation of $\pxw$ to a related polynomial $\pxtwt$.  We define
$\xti$ and $\wti$ now:
\begin{dfn}\label{dfn:tilde}
  Let
  \begin{equation*}
    \delxw = \{i: x(i) \neq w(i) \text{ or } \dxw(i,x(i)) \neq 0\}.
  \end{equation*}
  For $\delxw = \{d_1,d_2,\dotcoms d_k\}$ with $d_i < d_j$ for $i <
  j$, set
  \begin{align*}
    \xti &= \fl([x(d_1),x(d_2),\dotcoms x(d_k)]) \ \text{ and }\\
    \wti &= \fl([w(d_1),w(d_2),\dotcoms w(d_k)]).
  \end{align*}
\end{dfn}
Note that $\xti$ and $\wti$ are permutations in $\fsp{k}$.

\subsection{KL Polynomials}
\label{sec:klp}

While Kazhdan and Lusztig define the KL polynomials for general
Coxeter groups via the associated Hecke algebra, there is a purely
combinatorial description which we now give for the case of $\fsn$.
In order to give this definition succinctly, we let $[q^k]\pxw$ denote
the coefficient of $q^k$ in the polynomial $\pxw$, set
\begin{equation*}
  \mu(x,w) = [q^{(l(w)-l(x)-1)/2}]\pxw,
\end{equation*}
and define $c_s(x) = 1$ if $xs < x$; $c_s(x) = 0$ if $xs > x$. 
\begin{thm}[\cite{K-L1}] 
  There is a unique set of polynomials $\{\pxw\}_{x,w\in\fsn}$ such
    that, for all $x,w\in\fsn$:
  \begin{enumerate}
  \item $\pww = 1$
  \item $\pxw = 0$ when $x\not\leq w$.
  \item $\deg(P_{x,w}) \leq (l(w)-l(x)-1)/2$ when $x < w$
  \item For $s\in \simref$ with $ws < w$,
    \begin{equation*}
      \pxw = q^{c_s(x)}\pxws + q^{1-c_s(x)}\pxsws - \sumsb{z \leq
        ws\\ zs < z} \mu(z,ws)q^{\frac{l(w)-l(z)}{2}} \pxz.      
    \end{equation*}
  \end{enumerate}
\end{thm}
The proof, while not difficult, is intricate and we refer the
interested reader to the original paper \cite{K-L1} of Kazhdan and
Lusztig or to the more detailed exposition in Humphreys \cite{Hum}.
We note that $\mu(x,w)$ is the coefficient of the highest possible
power of $q$ in $\pxw$.

The KL polynomials satisfy many properties that are not immediately
apparent from the definition.  We list these properties without proof
below.  Properties 1 and 2 are standard and are due to Kazhdan and
Lusztig \cite{K-L1}.  (The polynomials $\pkl{w_0z}{w_0x}$ are referred
to as ``inverse'' KL polynomials.)  Property 3 is due to Lakshmibai
and Sandhya \cite{Lak-San} using results of Carrell \cite{carrell94}.
Properties 4 and 5 can be found in \cite{gwsb-msl}.

\begin{thm}\label{thm:stdfacts}
\begin{enumerate}\label{en:klprops}
\item For $s\in \simref$, $\pxw = \pxsw$ (resp., $\pxw = \psxw$) if
  $ws < w$ (resp., $sw < w$).\label{fact:pxwsim}
\item $\sum_{x\leq z\leq w} (-1)^{l(z)+l(w)} \pzw\pkl{w_0z}{w_0x} =
  \delta_{x,w}$ ($\delta$ is the Kronecker delta).\label{fact:invmat}
\item $\pxw = 1$ for all $x\leq w$ if and only if $w$ is 3412- and
  4231-avoiding.\label{fact:lsbeau}
\item $P_{\xti,\wti} = P_{x,w}$.\label{fact:pxw}
\item $\deg (\pxw) = 0$ if there do not exist $i < j < k < l$ such that 
  \begin{equation*}
    (\fl([x(i),x(j),x(k),x(l)]),\fl([w(i),w(j),w(k),w(l)]))
  \end{equation*}
  has a Bruhat picture of one of the two forms given in \figref{fig:2bsing}.
  \label{fact:2bsing}   
\end{enumerate}
\end{thm}
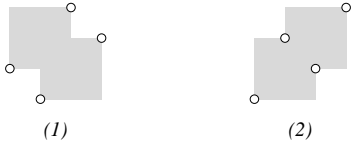
\mysmmslfig{2bsing}{Requirements for $\deg(\pxw) > 0$.}

We also have the following formulae for KL polynomials of certain
irreducible components of the singular loci of Schubert varieties.
The formulae are due to the author and Billey \cite{gwsb-msl} and,
independently, to both Manivel \cite{manivel2} and Cortez
\cite{cortez}.
\begin{thm}\label{thm:regkl}\ \\
  \begin{enumerate}
  \item For $k,m \geq 1$, 
    \begin{equation*}
      P_{\xkm,\wkm} = 1 + q + \cdots + q^{\min(k-1,m-1)}\label{eq:2}
  \end{equation*}
  and $\pkl{z}{\wkm} = 1$ for $\xkm < z \leq \wkm$.
  \item For $k,m \geq 1$,
    \begin{equation*}
      P_{\xktm,\wktm} = 1 + q \label{eq:6}
    \end{equation*}
    and
    $\pkl{z}{\wktm} = 1$ for $\xktm < z \leq \wktm$.
  \end{enumerate}
\end{thm}

The main ingredients in the proof of \thmref{thm:mainklinv} are
Theorems \ref{thm:regkl} and \ref{thm:stdfacts}.\ref{fact:invmat}.

\section[KL polynomials $\pkl{w_0w}{w_0x}$]
{KL polynomials $\pkl{w_0w}{w_0x}$}
\label{sec:kl-polyn-inverse}

Suppose $x\leq w$ such that $\pzw = 1$ for all $x < z \leq w$.  In
this case, Property \ref{fact:invmat} of \thmref{thm:stdfacts}
simplifies to
\begin{equation}\label{eq:1}
  \pkl{w_0w}{w_0x} = (-1)^{l(x) + l(w) + 1}\pxw +
  \sum_{x < z < w} (-1)^{l(z)+l(w)+1}\pkl{w_0z}{w_0x}.  
\end{equation}

Before utilizing \eqref{eq:1} to prove \thmref{thm:mainklinv}, we
first prove two technical lemmas.

\subsection{Two technical lemmas}
\label{sec:two-technical-lemmas}
\begin{lem}\label{lem:tech1} For $k,m \geq 1$,
  \begin{multline}\label{eq:4}
    \sum_{a,b=0}^{\minsb{a=k-1\\b=m-1}} (-1)^{a+b+1}
    \binom{k}{a}\binom{m}{b}
    \sum_{r=0}^{\minsb{\min(k-a-1,\\m-b-1)}}
    \binom{k-a-1}{r}\binom{m-b-1}{r} q^r\\ =
    (-1)^{k+m+1}\sum_{r=0}^{\minsb{\min(k-1,\\m-1)}}q^r.
  \end{multline}
\end{lem}
\begin{proof}
  It is convenient to sum $r$ from $0$ to $\min(k-1,m-1)$ rather than
  from $0$ to $\min(k-a-1,m-b-1)$.  As $\binom{n}{d} = 0$ whenever $0
  \leq n < d$, an extension of our summation range in this manner adds
  only terms equal to $0$.  So, we rewrite the left hand side of
  (\ref{eq:4}) as
  \begin{equation}\label{eq:5}
    \sum_{r=0}^{\minsb{\min(k-1,\\m-1)}} \sum_{a,b = 0}^{\minsb{a=k-1\\b=m-1}}
    (-1)^{a+b+1} \binom{k}{a}\binom{m}{b}
    \binom{k-a-1}{r}\binom{m-b-1}{r} q^r.
  \end{equation}
  It is a standard fact (see, e.g., \cite[(5.25)]{concmath}) that
  for $r \leq m-1$,
  \begin{equation*}
    \sum_{b=0}^{m-1} (-1)^b\binom{m}{b}\binom{m-b-1}{r} =
    (-1)^{r+m-1}\binom{m-r-1}{m-r-1} = (-1)^{r+m-1}.
  \end{equation*}
  Applying this identity twice to \eqref{eq:5}, we find that
  it equals
  \begin{equation*}
    \sum_{r=0}^{\min(k-1,m-1)} q^r (-1)^{r+k-1+1}(-1)^{r+m-1} = 
    \sum_{r=0}^{\min(k-1,m-1)} q^r (-1)^{k + m + 1}.
  \end{equation*}
\end{proof}
  
Our second lemma is similar in quality to the first, but requires an
additional piece of notation.  Define
\begin{equation*}
  f_{k,m}(a,b) = \binom{k}{a}\binom{m}{b}\left[(-1)^{a+b+1}\left(1
      + (k+m-a-b-1)q\right) + 2(-1)^{a+b}\right].
\end{equation*}

\begin{lem}\label{lem:tech2}
  For $k,m \geq 1$,
  \begin{equation}\label{eq:7}
    \sum_{a=0}^{k-1}\sum_{b=0}^{m-1} f_{k,m}(a,b) =
    (-1)^{k+m}(1+q).
  \end{equation}
\end{lem}
\begin{proof}
  We can rewrite the left side of (\ref{eq:7}) as
  \begin{equation}
    \label{eq:9}
    \sum_{a=0}^k \sum_{b=0}^m f_{k,m}(a,b) - \sum_{a=0}^k
    f_{k,m}(a,m) - \sum_{b=0}^m f_{k,m}(k,b) + f_{k,m}(k,m).
  \end{equation}
  By definition, $f_{k,m}(k,m) = (-1)^{k+m}(1+q)$.  Repeated
  application of the identities $\sum_{i=0}^n \binom{n}{i}(-1)^i = 0$
  and $\sum_{i=0}^n i \binom{n}{i}(-1)^i = 0$ (valid for $n > 0$)
  shows that the remaining terms in (\ref{eq:9}) are zero.
\end{proof}

\subsection{Inverse KL polynomials}
\label{sec:inverse-kl-polyn}
\begin{proof}[Proof of \thmref{thm:mainklinv}]
  As it streamlines the argument, we will prove the theorem for any
  pair $(x,w)$ for which $(\xti,\wti)$ equals $(\xkm,\wkm)$ or
  $(\ykm,\vkm)$ for some $k,m$.

  \noindent \textbf{Part 1.} We will argue by double induction on $k$
  and $m$.  First suppose that $(\xti,\wti) = (\xkm,\wkm)$ where
  either $k = 1$ or $m = 1$.  As can be seen by extrapolating from
  \figref{fig:pxw4231}.1, $\wkm$ is 3412- and 4231-avoiding in these
  cases.  Hence, by \thmref{thm:stdfacts}, parts \ref{fact:lsbeau} and
  \ref{fact:pxw}, $\pkl{\xkm}{\wkm} = 1$.  Now assume we have proven
  Part \ref{inveq:4231} for all pairs $x',w'$ for which $(\xti',\wti')
  = (\xrs,\wrs)$ with $1 \leq r \leq k$, $1 \leq s\leq m$ and $r+s <
  k+m$.
    
  Consider now the case of $x,w$ with $(\xti,\wti) = (\xkm,\wkm)$.
  There are two obvious simplifications we can make in the first term
  of (\ref{eq:1}).  First, we know from
  \thmref{thm:regkl}.\eqref{eq:2} that $\pxw = 1 + q + \cdots +
  q^{\min(k-1,m-1)}$.  Second, one can check that
  $l(w)-l(x) = k+m-1$.
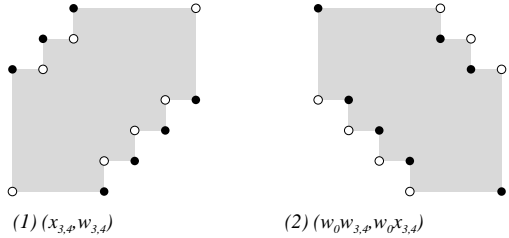
  \mysmmslfig{pxw4231}{We illustrate the action of multiplying a pair
    $(x,w)$ on the left by $w_0$.  Note that it corresponds to
    flipping the permutation matrix about a vertical line.}
  For the sum in (\ref{eq:1}), we start by noting that for any $z$
  with $x < z < w$, we must have $(\xti,\zti) = (\xrs,\wrs)$ for some
  $1\leq r\leq k$, $1\leq s\leq m$ and $r + s < k + m$ (see
  \figref{fig:pxw4231}).  We see that such a $z$ can be chosen in
  $\sum_{a,b}\binom{k}{a} \binom{m}{b}$ different ways.  Fix $r,s$ and
  set $a = k - r$ and $b = m - s$.  It follows that $a + b =
  l(w)-l(z)$.  This lets us rewrite $(-1)^{l(z)+l(w)+1}$ as $(-1)^{a +
  b + 1}$.
    
  Utilizing the above two facts and the induction hypothesis, we
  arrive at
  \begin{multline}\label{eq:3}
    \pkl{w_0w}{w_0x} = (-1)^{k+m}\sum_{r=0}^{\min(k-1,m-1)}q^r +
    \sum_{r=0} \binom{k-1}{r}\binom{m-1}{r} q^r +\\
    \sum_{a,b = 0}^{\minsb{a=k-1\\b=m-1}} (-1)^{a+b+1}
    \binom{k}{a}\binom{m}{b} \sum_{r=0}^{\min(k-a-1,m-b-1)}
    \binom{k-a-1}{r}\binom{m-b-1}{r} q^r.
  \end{multline}
  Note that the second sum in (\ref{eq:3}) is included to adjust for
  the fact that we allow $a = b = 0$ in the third term (i.e., the case
  of $z = w$).  By \lemref{lem:tech1}, the first and third terms
  cancel.  This proves Part 1.
    
  \bigskip
  \noindent
  \textbf{Part 2.} We will prove by double induction on $k$ and $m$.
  Define $y_{0,m} = x_{2,m}$, $y_{k,0} = x_{k,2}$, $v_{0,m} = w_{2,m}$
  and $v_{k,0} = w_{k,2}$.  The base cases of $k = 0$ or $m = 0$ (not
  both $0$) reduce to Part 1 of this theorem.  (In the ensuing
  induction, we do not use the case $k=m=0$.)  So for the remainder of
  the proof we assume $k,m \geq 1$.  Now assume we have proven
  (\ref{inveq:3412}) for all pairs $y',v'$ for which $(\yti',\vti') =
  (\xrts,\wrts)$ with $r \leq k$, $s\leq m$ and $0 < r+s < k+m$.
    
  Consider now the case of $y,v$ with $(\yti,\vti) = (\xktm,\wktm)$.
  Again, we will simplify (\ref{eq:1}).  We know by \thmref{thm:regkl}
  that $\pyv = 1 + q$.  Also, one can check that
  $l(v)-l(y) = k+m+1$.  These two facts let us write the first term in
  (\ref{eq:1}) as $(-1)^{k+m}(1+q)$.
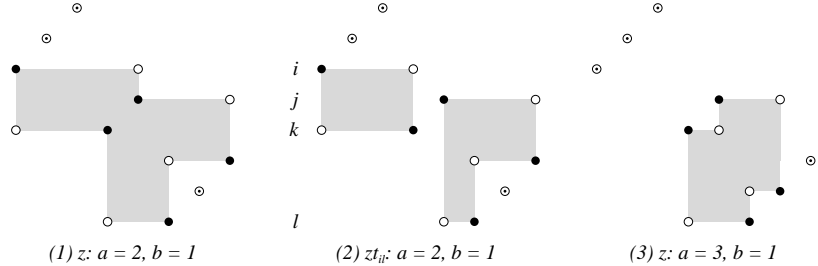
  \mysmmslfig{inv}{These three images illustrate for $k=m=3$ the first
    three types of $z$ terms arising in the proof of
    \thmref{thm:mainklinv}.2.  ($y$ is depicted with black dots and
    $z$ with white dots.)}
  We now categorize the $z$ for which $y < z < v$ and examine how they
  contribute to (\ref{eq:1}).  Note that by \thmref{thm:regkl}, $\pzv =
  1$ for each of these $z$.  Such $z$ fall into four categories.  See
  \figref{fig:inv} for examples of the different types.
  \begin{enumerate}
    
  \item For any $0 \leq a < k$ and $0 \leq b < m$, there are
    $\binom{k}{a}\binom{m}{b}$ different permutations $z$ with $y < z
    \leq v$ and $(\yti,\zti) = (y_{k-a,m-b},v_{k-a,m-b})$.  We know by
    the induction hypothesis that $\pkl{w_0z}{w_0y} = 1 +
    (k+m-a-b-1)q$ for such $z$.  Furthermore, it is easily checked
    that $(-1)^{l(z) + l(v)} = (-1)^{a + b}$.  We can write the
    contribution of these permutations $z$ succinctly in terms of
    $f_{k,m}(a,b)$.  However, as we are only interested in those $z <
    v$, we add in a corrective term to account for the fact that we
    allow $a = b = 0$ in our sum:

    \begin{equation}\label{eq:type1}
    1 + (k+m-1)q + \sum_{a,b=0}^{\minsb{a=k-1\\b=m-1}}
    \binom{k}{a}\binom{m}{b}\left[(-1)^{a+b+1}(1 + (k+m-a-b-1)q)\right].
    \end{equation}

  \item For each $z$ as in the previous case, $z\til$ and $z\tjk$ also
    lie strictly below $v$ and above $y$ (see \figref{fig:inv}).  By
    \thmref{thm:stdfacts}, parts \ref{fact:lsbeau} and
    \ref{fact:2bsing}, we know that $\pkl{w_0z\til}{w_0y} =
    \pkl{w_0z\tjk}{w_0y} = 1$.  Since $l(z\til) + l(w) = l(z\tjk) +
    l(v) = l(z) + l(v) + 1$, we get a total contribution of:

    \begin{equation}\label{eq:type2}
    \sum_{a,b=0}^{\minsb{a=k-1\\b=m-1}} 2(-1)^{a+b+1}.
    \end{equation}

    Note that we don't correct for $a = b = 0$ as $v\til, v\tjk < v$.
    
  \item There are $\binom{k}{a}\binom{m}{b}$ permutations
    $z$ for which $(\yti,\zti) = (y_{2-a,m-b},v_{2-a,m-b})$ for some
    $a$ and $b$ with $0\leq a \leq 1$ and $0\leq b < m$.  The
    cumulative contribution of these permutations $z$ to the sum in
    (\ref{eq:9}) can be determined using
    \thmref{thm:stdfacts}.\ref{fact:invmat}.  In particular,
    permutations of this form are precisely those lying in the
    interval $y < z \leq v'$ where $v' = wt_{12}t_{23}\cdots t_{k-1,k}
    t_{k,k+2}$.  We know from \thmref{thm:stdfacts}.\ref{fact:invmat}
    that
    \begin{equation*}
      \sum_{y \leq z \leq v'} (-1)^{l(z)} \pkl{z}{v'}
      \pkl{w_0z}{w_0y} = 0.
    \end{equation*}
    By \thmref{thm:stdfacts}.\ref{fact:lsbeau}, $\pkl{z}{v'} = 1$ for
    $y < z \leq v'$ and, by the first part of this theorem,
    $\pkl{y}{v'} = 1 + q$.  Hence, we see that these permutations
    contribute
    \begin{align*}
      (-1)^{l(v)+1}\sum_{y < z \leq v'} (-1)^{l(z)} \pkl{w_0z}{w_0y} &=
      (-1)^{l(v)+1}(-1)^{l(y)+1}(1+q)\\ &= (-1)^{k+m+1}(1+q).
    \end{align*}

  \item Similarly, the permutations $z$ for which $(\yti,\zti) =
    (y_{k-a,2-b},v_{k-a,2-b})$ also cumulatively contribute
    $(-1)^{k+m+1}(1+q)$.
  \end{enumerate}
  Incorporating the above knowledge into (\ref{eq:1}), we
  see that
  \begin{equation*}
    \begin{split}
      \pkl{w_0v}{w_0y} &= (-1)^{k+m}\pyv + \text{ (Type 1) } +
      \text{ (Type 2) } + \text{ (Type 3) } + \text{ (Type 4) }\\
      &= (-1)^{k+m}(1+q) + \\ &\left(1 + (k+m-1)q +
        \sum_{a,b=0}^{\minsb{a=k-1\\b=m-1}} f_{k,m}(a,b)\right) +\\
      &\ \ (-1)^{k+m+1}(1 + q) + (-1)^{k+m+1}(1 + q)\\
      &= 1 + (k+m-1)q + (-1)^{k+m+1}(1 + q) +
      \sum_{a=0}^{k-1}\sum_{b=0}^{m-1} f_{k,m}(a,b)\\
      &= 1 + (k+m-1)q.
    \end{split}
  \end{equation*}
\end{proof}

\bibliography{gen}
\end{document}